\RequirePackage{fix-cm}
\documentclass[smallextended]{svjour3}       
\smartqed  
\usepackage{amsmath,amssymb}
\usepackage{graphicx}
\usepackage{mathptmx}      
\begin{document}

\title{Modified hybrid combination synchronization of chaotic fractional order systems
}


\author{Kayode S Ojo         \and
        Samuel T Ogunjo        \and
        Ibiyinka A Fuwape
}


\institute{K. S. Ojo, S. T. Ogunjo, I. A. Fuwape \at
           Department of Physics \\
           Federal University of Technology, Akure\\
           Ondo State, Nigeria\\
           \email{stogunjo@futa.edu.ng}           
           \emph{Present address: Michael and Cecilia Ibru University, Ughelli, Delta State} of I. A. Fuwape  
}

\date{Received: date / Accepted: date}

\maketitle

\begin{abstract}
The paper investigates a new hybrid synchronization called modified hybrid synchronization (MHS) via the active control technique. Using the active control technique, stable controllers which enable the realization of the coexistence of complete synchronization, anti-synchronization and project synchronization in four identical fractional order chaotic systems were derived.  Numerical simulations were presented to confirm the effectiveness of the analytical technique.

\end{abstract}
\newpage
\section{Introduction}\label{intro}
A chaotic system is one whose motion is sensitive to initial conditions \cite{strogatz1994nonlinear}.  Since different initial conditions lead to different trajectories for the same dynamical system, it is expected that the trajectories cannot coincide. The possibility of two chaotic systems with different trajectories to follow the same trajectory by the introduction of a control function, as proposed by \cite{Pecora1990}, has been an interesting research area for scientists in nonlinear dynamics.  This is partly due to the applicability in different fields such as communication technology, security, neuroscience, atmospheric physics and electronics.

There are several methods for the synchronization of chaotic systems.  These methods include active control, Open Plus Closed Loop (OPCL), backstepping, feedback control, adaptive control, sliding mode and others.  A comparison of performance of a modified active control method and backstepping control  on synchronization of integer order system has been investigated \cite{ojo2013comparison}.  The active control method was found ``to be simpler with more stable synchronization time and hence more suitable for practical implementation".  The active control method was also found to have the best stability and convergence when compared with the direct method and OPCL method for fractional order systems \cite{ogunjo2017comparison}.

Generally, complete synchronization between a drive system $y_i$ and response system $x_i$ is said to occur if $\lim_{t\rightarrow+\infty}||y_i - x_i|| = 0$ and anti-synchronization if $\lim_{t\rightarrow+\infty}||y_i + x_i|| = 0$. If the error term is such that $\lim_{t\rightarrow+\infty}||y_i - \alpha x_i|| = 0$, where $\alpha$ is a positive integer, we have projective synchronization.  According to \cite{Zhou2017}, $\delta$ synchronization is defined by the error given as $\lim_{t\rightarrow+\infty}||y_i \pm x_i|| \leq \delta$, where $\delta$ has small value.  Other forms of synchronization include phase synchronization, anticipated synchronization, lag synchronization etc.  The possibility of one or more of these synchronization scheme in a single synchronization has not been explored.

Different synchronization methods and techniques have been used to study synchronization between two similar integer order systems \cite{ojo2012synchronization}, two dissimilar integer order systems of same dimension \cite{MOTALLEBZADEH20123643,femat2002}, two similar or dissimilar systems with different dimensions \cite{ojo2013mixed,ojo2014increased,ogunjo2013increased}, three or more integer order system (compound, combination-combination synchronization) \cite{ojo2014reduced1,ojo2014reduced2,Mahmoud2016}, discrete systems \cite{liu2008,kloeden2004,ma2007}, fractional order system of similar dimension \cite{LU20051125}, fractional order synchronization of different dimension \cite{Bhalekar2014,ayub2017},  circuit implementation of synchronization \cite{adelakun2017dynamics} and synchronization between integer order and fractional order systems \cite{chen2012}.

The study of chaotic systems has evolved over time from integer order dynamical systems to cover partial differential equations, time delayed differential equations, fractional order differential equations and even time series data.  The prevalence of integer order system was the lack of solution methods for fractional differential equations \cite{chen2009} and its inherent complexity \cite{gutierrez2010fractional}.  The Gr$\ddot{u}$nwald-Letnikov definition of fractional order systems, the fractional order derivative of order $\alpha$ can be written as \cite{petras2011fractional}
\begin{equation}\label{gl}
    D_t^\alpha f(t) = \lim_{h\rightarrow 0}\frac{1}{h^\alpha} \sum_{j=0}^\infty (-1)^j\begin{pmatrix}
                                                   \alpha \\
                                                   j\\
                                                 \end{pmatrix} f(t-jh)
\end{equation}
where the binomial coefficients can be written in terms of the Gamma function as
\begin{equation*}
    \begin{pmatrix}
    \alpha \\
    j\\
    \end{pmatrix}
    =
    \frac{\Gamma(\alpha + 1)}{\Gamma(j+1)\Gamma(\alpha - j +1)}
\end{equation*}
The Riemann-Liouville definition of fractional derivative is given as
\begin{equation}\label{rl}
    D_t^{-\alpha}f(t) = \frac{1}{\Gamma(n-\alpha)}  \frac{d^n}{dt^n} \int_a^t \frac{f(\tau)}{(t-\tau)^{\alpha +1}}d\tau
\end{equation}
The Caputo fractiional derivatives can be written as
\begin{equation}\label{caputo}
    D_t^{\alpha} f(t) = \frac{1}{\Gamma(n-\alpha)}\int_a^t \frac{f^{(n)}(\tau)}{(t-\tau)^{\alpha - n + 1}}d\tau,\,\, n-1<\alpha<n
\end{equation}
Fractional order systems have been found as a useful model in many engineering, physical and biological systems.

In this present work, we aim to investigate the possibility of coexistence of different synchronization scheme in the synchronization of four chaotic systems (two drives and two response systems).  Specifically, we aim to implement synchronization, anti-synchronization and projective synchronization on different dimensions in a fractional order combination synchronization using the method of active control.  We believe, if implemented, it will will enhance faster, robust and more secure information transmission.  To the best of our knowledge, this has not been reported in literature.
\section{System Description}\label{system}
The integer order Chen system was introduced by \cite{GuanrongChen1999} as
\begin{equation}\label{chen1}
    \begin{split}
      \dot{x} &= a(y - x) \\
      \dot{y} &= (c-a)x + cy - xz\\
      \dot{z} &= -bz + xy
    \end{split}
\end{equation}

\begin{figure}
  \includegraphics[width=\textwidth]{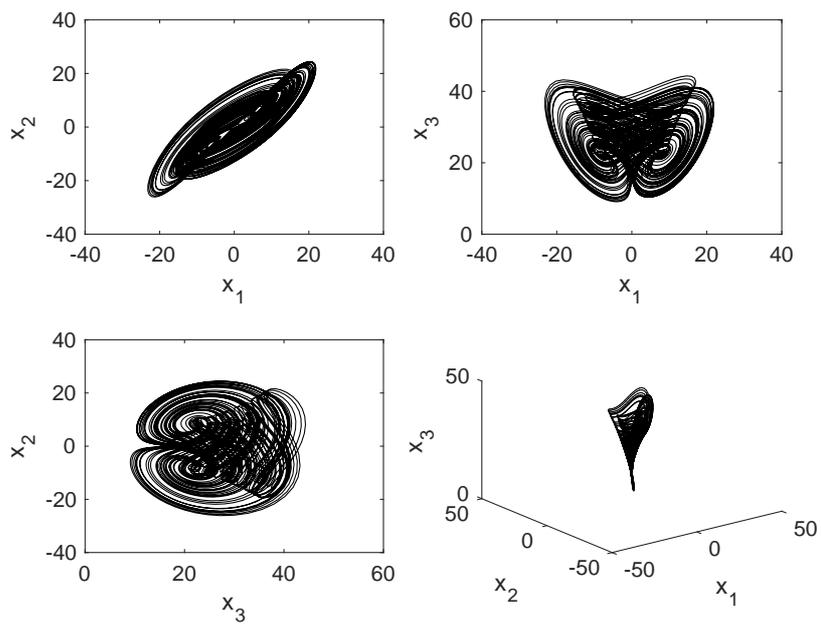}\\
  \caption{Phase space of the system. }\label{fig1}
\end{figure}

The fractional order chaotic Chen system was introduced by \cite{Li2004} as
\begin{equation}\label{eqn_d1}
    \begin{split}
      D^{\alpha} x_1 &= \sigma (x_2 - x_1) \\
      D^{\alpha} x_2 &= (c-a) x_1 - x_1x_3 + cx_2\\
      D^{\alpha} x_3 &= x_1 x_2 - b x_3
    \end{split}
\end{equation}
The system was found to be chaotic when $(a,b,c)=(35,3,28)$ and $0.7 \leq \alpha \leq 0.9$.  However, by varying parameter $a$ rather than parameter $c$ as in \cite{Li2004}, the system was found to be chaotic in the region $0.1\leq \alpha \leq 0.1$ \cite{LU2006685}. The phase space of the fractional order Chen system is shown in figure \ref{fig1}.   Various successful attempts have been made at synchronization of the integer, hyperchaotic, and fractional order Chen system \cite{hegazi2011,long2010,deng2005,Chen2008}.

\section{Design and implementation of synchronization scheme}\label{sec_des_imp}
The co-existence of different synchronization scheme within commensurate fractional order Chen system will be studied.  Suitable controllers are designed (Section \ref{design}) and numerical simulations presented in Section \ref{simulation} to verify the proposed controllers. \cite{Bhalekar2014}

\subsection{Design of controllers}\label{design}

Let the two drive system be defined as
\begin{equation}\label{eqn_d1}
    \begin{split}
      D^{p_1} x_1 &= \sigma (x_2 - x_1) \\
      D^{p_2} x_2 &= (c-a) x_1 - x_1x_3 + cx_2\\
      D^{p_3} x_3 &= x_1 x_2 - b x_3
    \end{split}
\end{equation}
and
\begin{equation}\label{eqn_d2}
    \begin{split}
      D^{q_1} y_1 &= \sigma (y_2 - y_1) \\
      D^{q_2} y_2 &= (c-a) y_1 - y_1y_3 + cy_2\\
      D^{q_3} y_3 &= y_1 y_2 - b y_3
    \end{split}
\end{equation}
Defining the two response systems as
\begin{equation}\label{eqn_d3}
    \begin{split}
      D^{r_1} z_1 &= \sigma (z_2 - z_1) + u_1 \\
      D^{r_2} z_2 &= (c-a) z_1 - z_1z_3 + cz_2 + u_2\\
      D^{r_3} z_3 &= z_1 z_2 - b z_3 + u_3
    \end{split}
\end{equation}
and
\begin{equation}\label{eqn_d4}
    \begin{split}
      D^{s_1} w_1 &= \sigma (w_2 - w_1) + u_4 \\
      D^{s_2} w_2 &= (c-a) w_1 - w_1w_3 + cw_2 + u_5\\
      D^{s_3} w_3 &= w_1 w_2 - b w_3 + u_6
    \end{split}
\end{equation}
where the six active control functions $u_1,\, u_2,\, u_3,\, u_4,\, u_5,\, u_6$ introduced in equations \ref{eqn_d3} and \ref{eqn_d4} are control functions to be determined.

We define the error states $e_1,\, e_2\, e_3$ as

\begin{equation}\label{eqn_d5}
    \begin{split}
      e_1 &= (x_1 + y_1) - (z_1 + w_1) \\
     e_2 &= (x_2 + y_2) + (z_2 + w_2) \\
     e_3 &= (x_3 + y_3) - \alpha (z_3 + w_3)
    \end{split}
\end{equation}
Substituting the drive systems (equations \ref{eqn_d1} and \ref{eqn_d2}) and response systems (equations \ref{eqn_d3} and \ref{eqn_d4}) into equation \ref{eqn_d5} and assuming a commensurate system, the error system is obtained as

\begin{equation}\label{eqn_d6}
    \begin{split}
      D^\mu e_1 &= -[ ae_1 + ae_2 - 2a(x_2 + y_2) + u_1 + u_4]\\
      D^\mu e_2 &= (a+c)e_1 + ce_2 + 2c(z_1 + w_1) - 2a(x_1 + y_1) - x_1x_3 - y_1y_3 -z_1z_3 - w_1w_3 + u_2 + u_5 \\
      D^\mu e_3 &= -[be_3 - x_1x_2 - y_1y_2 + \alpha z_1z_2 + \alpha w_1w_2 + \alpha u_3 + \alpha u_6
    \end{split}
\end{equation}
Active control inputs $u_i (i=1,2,3,4,5,6)$ are then defined as
\begin{equation}\label{eqn_d7}
    \begin{split}
      u_1 + u_4 &=  -[V_1 + 2a(x_2 + y_2) \\
      u_2 + u_5 &= V_2 - 2c(z_1 + w_1) + 2a(x_1 + y_1) + x_1x_3 + y_1y_3 + z_1z_3 + w_1w_3\\
      u_3 + u_6 &= \frac{1}{\alpha}[-V_3 + x_1x_2 + y_1y_2 - \alpha z_1z_2 - \alpha w_1w_2 ]
    \end{split}
\end{equation}
where the functions $V_i$ are to be obtained.  Substituting equation \ref{eqn_d7} into equation \ref{eqn_d6} yields
\begin{equation}\label{eqn_d8}
    \begin{split}
      D^\mu e_1 &= -ae_1 - ae_2 + V_1\\
      D^\mu e_2 &= (a+c)e_1 + ce_2 + V_2 \\
      D^\mu e_3 &= -be_3 + V_3
    \end{split}
\end{equation}
The synchronization error system (equation \ref{eqn_d8}) is a linear system with active control inputs $V_i$.  We design  an appropriate feedback control which stabilizes the system so that $e_i (i = 1,2,3) \rightarrow 0$ as $t\rightarrow\infty$, which implies that synchronization is achieved with the proposed feedback control.  There are many possible choices for the control inputs $V_i$, for simplicity, we chose
\begin{equation}\label{eqn_d9}
  \begin{bmatrix}
    V_1 \\
    V_2 \\
    V_3 \\
  \end{bmatrix}
  = C \begin{bmatrix}
    e_1 \\
    e_2 \\
    e_3 \\
  \end{bmatrix}
\end{equation}
where C is a $3\times3$ constant matrix.  In order to make the closed loop system stable, matrix C should be selected in such a way that the feedback system has eigenvalues $\lambda_i$ that satisfies the equation
\begin{equation}\label{eqn_d10}
    |\arg(\lambda_i)| > 0.5\pi\alpha, \,\, i = 1,2,\ldots.
\end{equation}
where $\lambda$ is the eigenvalue, I is an identity matrix and A is the coeffcient of the error state.  There are varieties of choices for choosing matrix C.  Matrix C is chosen as

\begin{equation}\label{eqn_d11}
    C = \begin{pmatrix}
          (c-\lambda) & a & 0 \\
          -(a+c) & -(c+\lambda) & 0 \\
          0 & 0 & (b-\lambda) \\
        \end{pmatrix}
\end{equation}
Using equation \ref{eqn_d11} in \ref{eqn_d9}, we obtain our control function as
\begin{equation}\label{eqn_d12}
    \begin{split}
      u_1 + u_4 &= -V_1 + 2a(x_2 + y_2) \\
      u_2 + u_5 &= V_2 - 2c(z_1 + w_1) + 2a(x_1 + y_1) + x_1x_3 + y_1y_3 + z_1z_3 + w_1w_3\\
      u_3 + u_6 &= \frac{1}{\alpha}[-V_3 + x_1x_2 + y_1y_2 - \alpha z_1z_2 - \alpha w_1w_2]
    \end{split}
\end{equation}

Based on the controllers obtained, two unique cases can be observed.

The control system can be defined as
\begin{equation}\label{eqn_d13}
    \begin{split}
      u_1  &= \frac{1}{2}[-V_1 + 2a(x_2 + y_2)] \\
      u_2  &= \frac{1}{2}[V_2 - 2c(z_1 + w_1) + 2a(x_1 + y_1) + x_1x_3 + y_1y_3 + z_1z_3 + w_1w_3]\\
      u_3  &= \frac{1}{2\alpha}[-V_3 + x_1x_2 + y_1y_2 - \alpha z_1z_2 - \alpha w_1w_2]\\
      u_4  &= u_1\\
      u_5  &= u_2\\
      u_6  &= u_3
    \end{split}
\end{equation}
It can also be defined as
The control system can be defined as
\begin{equation}\label{eqn_d14}
    \begin{split}
      u_1  &= -V_1 + 2a(x_2 + y_2)\\
      u_2  &= V_2 - 2c(z_1 + w_1) + 2a(x_1 + y_1) + x_1x_3 + y_1y_3 + z_1z_3 + w_1w_3\\
      u_3  &= \frac{1}{\alpha}[-V_3 + x_1x_2 + y_1y_2 - \alpha z_1z_2 - \alpha w_1w_2]\\
      u_4  &= 0\\
      u_5  &= 0\\
      u_6  &= 0
    \end{split}
\end{equation}

\subsection{Numerical simulation of Results}\label{simulation}
\begin{figure}
  \includegraphics[width=\textwidth]{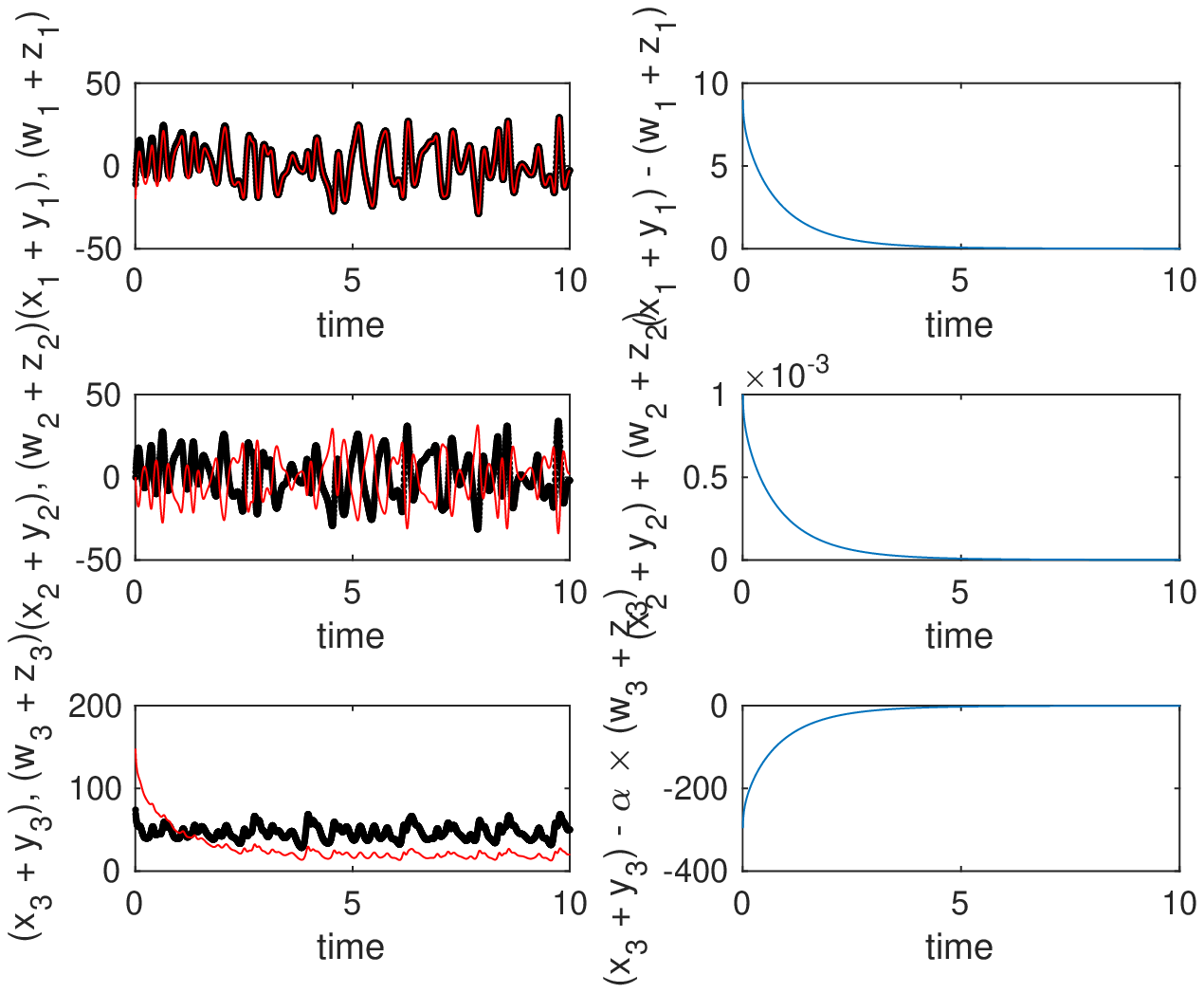}\\
  \caption{Synchronization obtained from the realization of case 1. }\label{fig2}
\end{figure}

\begin{figure}
  \includegraphics[width=\textwidth]{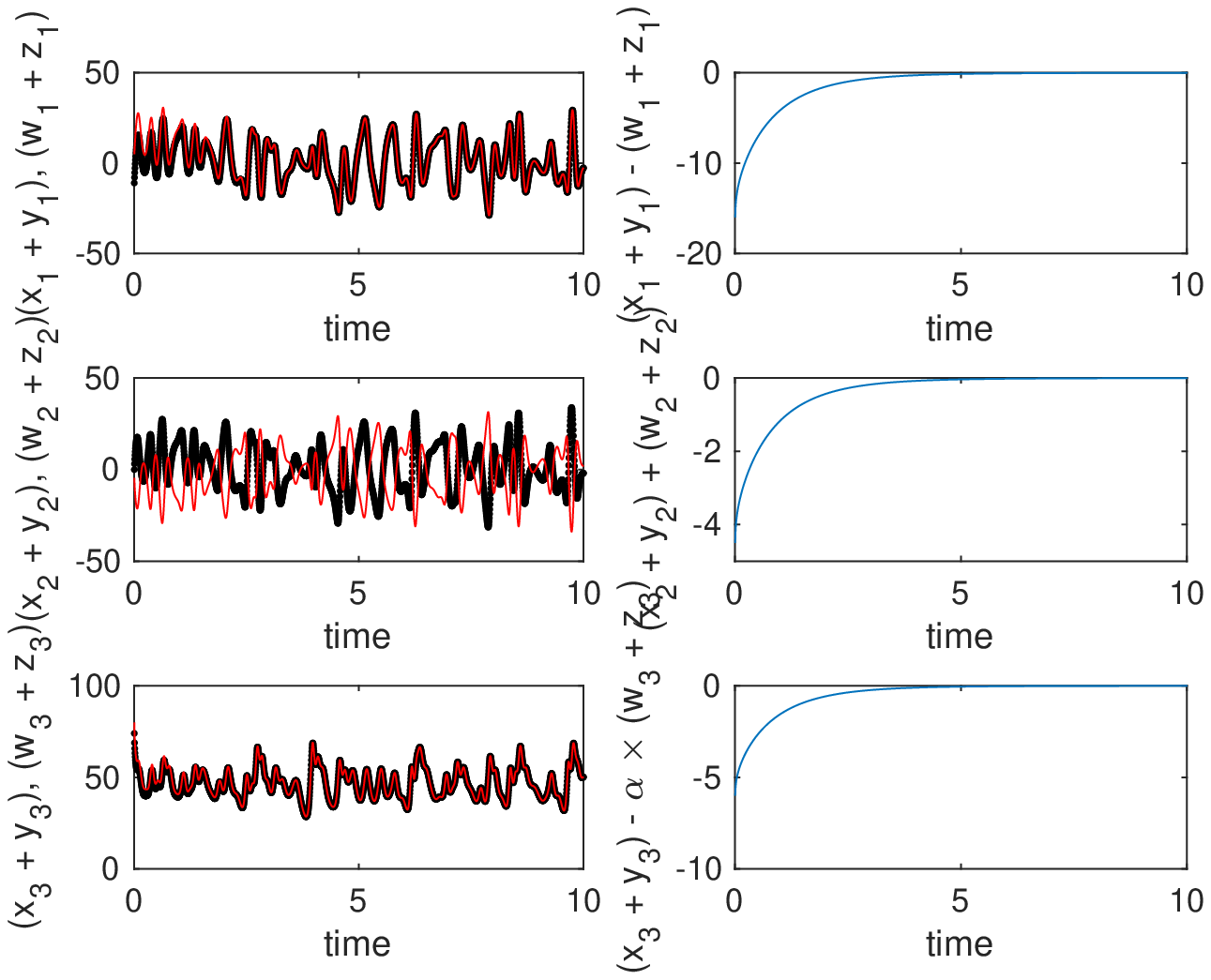}\\
  \caption{Synchronization obtained from the realization of case 2. }\label{fig3}
\end{figure}

To verify the effectiveness of the synchronization scheme proposed in section \ref{design} using the method of active control, we used the initial conditions $x_i(-10,0.001,37)$, $y_i(37,-5,0)$, $w_i(-5,0.5,25)$ and $z_i(10,-5,15)$.  The order of the system was taken as 0.95.  A time step of 0.005 was used.  In the case of projective synchronization, the scaling parameter was taken to be 5. The parameters of the system are taken as $(a,b,c) = (35,3,28)$.  According to \cite{petras2011fractional}, the general numerical solution of the fractional differential equation
\begin{equation}\label{eqn_ff}
    _aD^q_t y(t) = f(y(t),t)
\end{equation}
can be expressed as
\begin{equation}\label{eqn_fg}
    y(t_k) = f(y(t_k),t_k)h^q - \sum_{j=v}^k c_j^{(q)} y(t_{t-j})
\end{equation}
where $c_i^{(q)}$ is given as
\begin{equation}\label{eqn_fh}
    \begin{split}
      c_0^{(q)} &= 1 \\
      c_j^{(q)} &= \left( 1 - \frac{1+q}{j}c_{j-1}^q \right)
    \end{split}
\end{equation}
The results for the two cases considered are shown in figures \ref{fig2} and \ref{fig3}.   From the results presented, the drives and responses were found to achieve synchronization as indicated by the convergence of the error terms to zero.  The effectiveness of the proposed scheme is hereby confirmed.

\section{Conclusion}
In this paper, a new synchronization scheme is proposed and implemented.  The modified hybrid synchronization that allows for the coexistence of different synchronization schemes was implemented in a compound synchronization of fractional order Chen system.  In particular, the controllers consists of complete synchronization, anti-synchronization, and projective synchronization.  We believe that this type of synchronization will offer better security and more robust.  There is the need to investigate the performance of this type of synchronization using different synchronization schemes.  Furthermore, it will be productive to study the behaviour of this scheme under different types and strength of noise.  Practical implementation of this scheme is also proposed.

\section*{Conflict of Interest}
The authors hereby declare that there is no conflict of interest.

\clearpage
\newpage

\bibliographystyle{spmpsci}      
\bibliography{my_bib}

\end{document}